\documentclass[10pt]{article}
\usepackage{mystyle}
\usepackage{mathtools}
\usepackage{youngtab}
\usepackage{xypic}
\title{Quiver relations and associated symmetric polynomials}
\author{Dimitry Noshchenko}
\begin{document}
\maketitle

\begin{abstract}
The idea is to identify certain path algebra elements with symmetric functions. We propose such a morphism
by solving the quiver relations, which describe the Plucker-type embedding for quiver grassmannians.
\end{abstract}

\section{Introduction}
Quiver is a directed graph with multiply edges and loops allowed.
Let $Q$ = finite quiver with $k$ vertices and $l$ arrows.
Quiver representation $M$ is a collection of vector spaces and linear maps, attached to each vertex
and arrow of $Q$:
\begin{equation}
M=\{V_1,\dots,V_k\}\cup \{M_{v,i,j}\}, \quad \mathrm{dim}(V_i)=d_i,\quad M_{v,i,j}=Mat(d_i\times d_j,\mathrm{k}): V_i \rightarrow V_j
\end{equation}
Quiver grassmannian $Gr_{\underline{e}}(Q,M)$ of dimension vector $\underline{e}=(e_1,\dots,e_k)$ contains all
$\underline{e}$-dimensional subspaces of $M$, compatible with the underlying linear maps.

Fix an ordered basis $\mathcal{B}=\{1,2,\dots,d_1,d_1+1,\dots,d_l\}$, $I\subset B$ for all $V_i$s. Quiver relations (QR)  are
of the form \cite{4}:
\begin{equation}\label{QR}
\sum_{s',t'} (-1)^{\epsilon(s',t')}M_{v,s',t'} \ \Delta_s^{s'} \ \Delta_{t'}^t 
- \ \sum_{s'} (-1)^{\epsilon(s',t)} M_{v,s',t} \ \Delta_{s}^{s'} \ \Delta_I = 0, \quad I\subset\mathcal{B}, \quad |I|=e
\end{equation}
\textbf{Define $Fset$ as a collection of all non-vanishing $\Delta^i_j$'s in \ref{QR}.}
Our aim is to generate polynomials which are sums of $GL$ characters with coefficients given by the embedding
\begin{equation}\label{Emb}
Gr_{\overline(e)}(Q,M)\xrightarrow{\phi_1} Gr(e,d) \xrightarrow{\phi_2} \mathbb{P}^{e+d-1}, \quad e=\sum e_i, d=\sum d_i
\end{equation}
Define 2 kind of invariants:
\begin{equation}
P(Q,s,\mathbf{M})=\sum\limits_{R}\Delta_R\chi_R(\overline{t})
\end{equation}
with $\Delta=$ homogeneous coordinate on the image of embedding \ref{Emb}, and
\begin{equation}
\tilde{P}(Q,s,\mathbf{M})=\sum\limits_{R'}\tilde{\Delta}_{R'}\chi_R(\overline{t})
\end{equation}
with $\tilde{\Delta}\in Fset$ only (PR are not taken into part).
Note that $P,P'$ are inhomogeneous in $R,R'$ (may contain Young diagrams of different size).
\\
\textbf{Proposition 1.} \emph{In case of $\chi_{\lambda}=s_{\lambda}$ polynomials $P(Q,s,\mathbf{M})$ form a basis of Symm algebra.} This seems to be obvious, but not the only one possible choice of $\chi$.
\\
Write $\mathbf{M}$ as block matrix:
\begin{equation}
\mathbf{M}= \left[ \begin {array}{cccc} M_{{1,1}}&M_{{1,2}}&M_{{1,3}}&\dots
\\\noalign{\medskip}M_{{2,1}}&M_{{2,2}}&M_{{2,3}}&
\\\noalign{\medskip}M_{{3,1}}&M_{{3,2}}&M_{{3,3}}&
\\\noalign{\medskip}\vdots& & &\ddots\end {array} \right] 
\end{equation}
with each $M_{i,j}$ corresponding to $(i,j)$ linear map.

\section{Computation of $P,\tilde{P}$. Quiver relations as vanishing minors}

\subsection{$Q=A_2$, $M=\{\mathbb{C},\mathbb{C},\mu=Mat(2,2,\mathbb{C})\},\overline{e}=(1,1)$}
$I=\{2,4\}$:
One non-trivial equation,  $\#(+,-)=(2,2)$
\begin{equation}
F \left( v,2,3 \right)=\mu_{{v,1,4}}\Delta_{{ \left\{ 1,4 \right\} }}\Delta_{{ \left\{ 2,3
 \right\} }}+\mu_{{v,2,4}}\Delta_{{ \left\{ 2,4 \right\} }}\Delta_{{
 \left\{ 2,3 \right\} }}-\mu_{{v,1,3}}\Delta_{{ \left\{ 1,4 \right\} }
}\Delta_{{ \left\{ 2,4 \right\} }}-\mu_{{v,2,3}}{\Delta_{{ \left\{ 2,4
 \right\} }}}^{2}=0
\end{equation}
Or:
\begin{equation}
F \left( v,2,3 \right)=\left| \begin {array}{ccc} \mu_{{v,1,4}}&0&\Delta_{{ \left\{ 2,4
 \right\} }}\\\noalign{\medskip}\mu_{{v,2,4}}&\Delta_{{ \left\{ 1,4
 \right\} }}&{\frac {\mu_{{v,2,3}}{\Delta_{{ \left\{ 2,4 \right\} }}}^
{2}}{\mu_{{v,1,4}}\Delta_{{ \left\{ 2,3 \right\} }}}}
\\\noalign{\medskip}\mu_{{v,1,3}}&\Delta_{{ \left\{ 2,3 \right\} }}&
\Delta_{{ \left\{ 2,3 \right\} }}\end {array} \right|=
\left| \begin {array}{ccc} \mu_{{v,1,4}}&\Delta_{{ \left\{ 2,4
 \right\} }}&\Delta_{{ \left\{ 2,4 \right\} }}\\\noalign{\medskip}-\mu
_{{v,2,4}}&\Delta_{{ \left\{ 1,4 \right\} }}&-{\frac {\mu_{{v,2,3}}
\Delta_{{ \left\{ 2,4 \right\} }}}{\mu_{{v,1,3}}}}\\\noalign{\medskip}
\mu_{{v,1,3}}&0&\Delta_{{ \left\{ 2,3 \right\} }}\end {array} \right| 
 =0
\end{equation}
\begin{equation}
,\quad \Delta_{1,2}=0,
\quad \Delta_{3,4}=0
\end{equation}
Note that if we inverse the arrow, our equation will change:
\begin{equation}
F \left( v,4,1 \right) =\mu'  _{{v,3,2}}
\Delta_{{ \left\{ 2,3 \right\} }}\Delta_{{ \left\{ 1,4 \right\} }}+
\mu'  _{{v,4,2}}\Delta_{{ \left\{ 2,4
 \right\} }}\Delta_{{ \left\{ 1,4 \right\} }}-\mu'
  _{{v,3,1}}\Delta_{{ \left\{ 2,3 \right\} }}\Delta_{{
 \left\{ 2,4 \right\} }}-\mu'  _{{v,4,1}
}{\Delta_{{ \left\{ 2,4 \right\} }}}^{2}=0,
\end{equation}
where $\mu'$ is the inverse map.
Now fix $\mu= \left( \begin {array}{cc} p&1-p\\\noalign{\medskip}1-p&p\end {array}
 \right) $ -- quantum permutation matrix from https://arxiv.org/abs/1109.4888.
Combining this with the only Plu"cker identity 
\begin{equation}
\Delta_{{ \left\{ 1,2 \right\} }}\Delta_{{ \left\{ 3,4 \right\} }}-
\Delta_{{ \left\{ 1,3 \right\} }}\Delta_{{ \left\{ 2,4 \right\} }}+
\Delta_{{ \left\{ 1,4 \right\} }}\Delta_{{ \left\{ 2,3 \right\} }}=0
\end{equation}
we obtain our polynomial invariant:
\begin{equation}
\boxed{P(A_2,4,\mu)=Qh_{{1}}-{\frac {Q \left( -Wp-Q+pQ \right) h_{{2}}}{W-Wp+pQ}}+W
 \left( {h_{{1}}}^{2}-h_{{2}}h_{{0}} \right) -{\frac { \left( -Wp-Q+pQ
 \right) W \left( h_{{1}}h_{{2}}-h_{{3}}h_{{0}} \right) }{W-Wp+pQ}}}
\end{equation}
where $Q,W$ are arbitrary complex numbers.

\subsection{$M=\{\mathbb{C}^2,\mathbb{C}^2,\mu=Mat(4,4,\mathbb{C})\},\overline{e}=(2,2)$}
$I={3,4,7,8}$: We have 2x2 non-trivial quiver equations () with $\#(+,-)=(4,5),(6,3)$
One of them ()
\begin{equation}
F(v,3,5)=
 \left| \begin {array}{ccc} \mu_{{v,1,7}}+\mu_{{v,2,7}}+\mu_{{v,2,8}}&
\Delta_{{ \left\{ 3,4,7,8 \right\} }}&\Delta_{{ \left\{ 3,4,5,7
 \right\} }}+2\,\Delta_{{ \left\{ 3,4,7,8 \right\} }}
\\\noalign{\medskip}\mu_{{v,3,8}}-\mu_{{v,3,7}}-\mu_{{v,1,5}}&\Delta_{
{ \left\{ 1,4,7,8 \right\} }}+2\,\Delta_{{ \left\{ 2,4,7,8 \right\} }}
&0\\\noalign{\medskip}-\mu_{{v,1,8}}+\mu_{{v,2,5}}+{\frac {\mu_{{v,3,5
}}\Delta_{{ \left\{ 3,4,7,8 \right\} }}}{\Delta_{{ \left\{ 2,4,7,8
 \right\} }}}}&\Delta_{{ \left\{ 3,4,7,8 \right\} }}+\Delta_{{
 \left\{ 1,4,7,8 \right\} }}&2\,\Delta_{{ \left\{ 3,4,5,8 \right\} }}+
\Delta_{{ \left\{ 3,4,5,7 \right\} }}\end {array} \right|=0
\end{equation}
\[
\mu_{{v,2,8}}\Delta_{{ \left\{ 2,4,7,8 \right(\} }}\Delta_
{{ \left\{ 3,4,5,7 \right\} }}+\mu_{{v,3,8}}\Delta_{{ \left\{ 3,4,7,8
 \right\} }}\Delta_{{ \left\{ 3,4,5,7 \right\} }}-\mu_{{v,1,5}}\Delta_
{{ \left\{ 1,4,7,8 \right\} }}\Delta_{{ \left\{ 3,4,7,8 \right\} }}
\]
\[ -
\mu_{{v,2,5}}\Delta_{{ \left\{ 2,4,7,8 \right\} }}\Delta_{{ \left\{ 3,
4,7,8 \right\} }}-\mu_{{v,3,5}}{\Delta_{{ \left\{ 3,4,7,8 \right\} }}}
^{2}
\]
Can we rewrite them in a determinantal form? 

{cont.}
\begin{equation}
F_(v,4,5)= \left| \begin {array}{ccc} \mu_{{v,4,7}}&m_{2;{1,2}}&\Delta_
{{ \left\{ 3,4,5,7 \right\} }}\\\noalign{\medskip}0&\Delta_{{ \left\{ 
3,4,7,8 \right\} }}&0\\\noalign{\medskip}-\mu_{{v,4,8}}&m_{2;{3,2}}&
\Delta_{{ \left\{ 3,4,5,8 \right\} }}\end {array} \right| + \left| 
\begin {array}{ccc} \mu_{{v,2,5}}&0&\Delta_{{ \left\{ 3,4,5,8
 \right\} }}\\\noalign{\medskip}-\mu_{{v,1,7}}&\Delta_{{ \left\{ 2,3,7
,8 \right\} }}&0\\\noalign{\medskip}\mu_{{v,2,7}}&\Delta_{{ \left\{ 1,
3,7,8 \right\} }}&\Delta_{{ \left\{ 3,4,7,8 \right\} }}\end {array}
 \right| +
\end{equation}

\[
+ \left| \begin {array}{ccc} -\mu_{{v,1,8}}&\Delta_{{
 \left\{ 2,3,7,8 \right\} }}&\Delta_{{ \left\{ 3,4,7,8 \right\} }}
\\\noalign{\medskip}\mu_{{v,2,8}}&\Delta_{{ \left\{ 1,3,7,8 \right\} }
}&0\\\noalign{\medskip}-\mu_{{v,1,5}}&-{\frac {\mu_{{v,4,5}}\Delta_{{
 \left\{ 3,4,7,8 \right\} }}}{\mu_{{v,2,8}}}}&\Delta_{{ \left\{ 3,4,5,
7 \right\} }}\end {array} \right|=0.
\]

{Polynomial $P(A_2,4,Id)$}
define $P(Q, n, M) =$ linear combination of GL characters with coefficients running through the Fset.
\textbf{this turns out to be a special base of Symm.}
This case gives
\begin{equation}
Fset=\left( \begin {array}{cccccccc} &&\Delta_{{ \left\{ 1,4,7,8
 \right\} }}&\Delta_{{ \left\{ 1,3,7,8 \right\} }}&&&&
\\\noalign{\medskip}&&\Delta_{{ \left\{ 2,4,7,8 \right\} }}&\Delta_{
{ \left\{ 2,3,7,8 \right\} }}&&&&\\\noalign{\medskip}&&\Delta_{{
 \left\{ 3,4,7,8 \right\} }}&&&&&\\\noalign{\medskip}&&&\Delta
_{{ \left\{ 3,4,7,8 \right\} }}&&&&\\\noalign{\medskip}&&&&&
&\Delta_{{ \left\{ 1,4,7,8 \right\} }}&-\Delta_{{ \left\{ 1,3,7,8
 \right\} }}\\\noalign{\medskip}&&&&&&\Delta_{{ \left\{ 2,4,7,8
 \right\} }}&-\Delta_{{ \left\{ 2,3,7,8 \right\} }}
\\\noalign{\medskip}&&&&&&\Delta_{{ \left\{ 3,4,7,8 \right\} }}&
\\\noalign{\medskip}&&&&&&&\Delta_{{ \left\{ 3,4,7,8 \right\} }
}\end {array} \right)
\end{equation}
 
Now we can write down the explicit form of $P(A_2,4,Id)$:
\begin{equation}
{\it Fset}=[\Delta_{{ \left\{ 1,3,7,8 \right\} }},\Delta_{{ \left\{ 1,
4,7,8 \right\} }},\Delta_{{ \left\{ 2,3,7,8 \right\} }},\Delta_{{
 \left\{ 2,4,7,8 \right\} }},\Delta_{{ \left\{ 3,4,5,7 \right\} }},
\Delta_{{ \left\{ 3,4,5,8 \right\} }},\Delta_{{ \left\{ 3,4,6,7
 \right\} }},\Delta_{{ \left\{ 3,4,6,8 \right\} }},\Delta_{{ \left\{ 3
,4,7,8 \right\} }}]
\end{equation}
Assuming all free parameters $=\pm1$, we get the following solution:
\begin{equation}
{\it Fset}=[1, 1, 1, -1, -1, 1, -1, -1, 1]
\end{equation}

Now choose $\mu=\mu(p,q)$. Solving QR's independently, we get:
\begin{equation}
\Delta_{{ \left\{ 1,3,7,8 \right\} }}=-\frac{1}{pq}(\Delta_{{ \left\{ 3,4,5,7
 \right\} }}-2\,\Delta_{{ \left\{ 3,4,5,7 \right\} }}q-\Delta_{{
 \left\{ 1,4,7,8 \right\} }}p-\Delta_{{ \left\{ 2,4,7,8 \right\} }}+p
\Delta_{{ \left\{ 2,4,7,8 \right\} }}+
\end{equation}
\[
+\Delta_{{ \left\{ 1,4,7,8
 \right\} }}pq+\Delta_{{ \left\{ 2,4,7,8 \right\} }}q-p\Delta_{{
 \left\{ 2,4,7,8 \right\} }}q-\Delta_{{ \left\{ 2,3,7,8 \right\} }}q+p
\Delta_{{ \left\{ 2,3,7,8 \right\} }}q)
\]
\begin{equation}
\Delta_{{ \left\{ 3,4,5,8 \right\} }}=-\frac{1}{q}(\Delta_{{ \left\{ 3,4,5,7
 \right\} }}-\Delta_{{ \left\{ 3,4,5,7 \right\} }}q-\Delta_{{ \left\{ 
1,4,7,8 \right\} }}p-\Delta_{{ \left\{ 2,4,7,8 \right\} }}+p\Delta_{{
 \left\{ 2,4,7,8 \right\} }})
\end{equation}
\begin{equation}
\Delta_{{ \left\{ 3,4,6,7 \right\} }}=-\frac{1}{p(2q-1)}(2\,p\Delta_{{ \left\{ 2,3,7,8
 \right\} }}q+2\,\Delta_{{ \left\{ 3,4,5,7 \right\} }}pq-2\,p\Delta_{{
 \left\{ 2,4,7,8 \right\} }}q-p\Delta_{{ \left\{ 3,4,5,7 \right\} }}+2
\,p\Delta_{{ \left\{ 2,4,7,8 \right\} }}-
\end{equation}
\[
-2\,\Delta_{{ \left\{ 3,4,5,7
 \right\} }}q+\Delta_{{ \left\{ 2,4,7,8 \right\} }}q-\Delta_{{
 \left\{ 2,3,7,8 \right\} }}q+\Delta_{{ \left\{ 3,4,5,7 \right\} }}-
\Delta_{{ \left\{ 2,4,7,8 \right\} }})
\]
\begin{equation}
\Delta_{\{3, 4, 6, 8\}}=\frac{1}{pq(2q-1)}(\Delta_{{ \left\{ 3,4,6,8 \right\} }}=2\,p\Delta_{{ \left\{ 2,4,7,8
 \right\} }}-2\,\Delta_{{ \left\{ 3,4,5,7 \right\} }}{q}^{2}p-2\,
\Delta_{{ \left\{ 1,4,7,8 \right\} }}{p}^{2}q+2\,{p}^{2}\Delta_{{
 \left\{ 2,4,7,8 \right\} }}q+2\,\Delta_{{ \left\{ 2,4,7,8 \right\} }}
{q}^{2}p-
\end{equation}
\[
-
2\,\Delta_{{ \left\{ 2,3,7,8 \right\} }}{q}^{2}p-p\Delta_{{
 \left\{ 3,4,5,7 \right\} }}+
\Delta_{{ \left\{ 1,4,7,8 \right\} }}{p}^
{2}-{p}^{2}\Delta_{{ \left\{ 2,4,7,8 \right\} }}+3\,\Delta_{{ \left\{ 
3,4,5,7 \right\} }}pq-4\,p\Delta_{{ \left\{ 2,4,7,8 \right\} }}q+
\]
\[
+2\,
\Delta_{{ \left\{ 1,4,7,8 \right\} }}pq+2\,p\Delta_{{ \left\{ 2,3,7,8
 \right\} }}q+2\,\Delta_{{ \left\{ 2,4,7,8 \right\} }}q-\Delta_{{
 \left\{ 2,3,7,8 \right\} }}q+
 \]
 \[ 
 \Delta_{{ \left\{ 3,4,5,7 \right\} }}-
\Delta_{{ \left\{ 2,4,7,8 \right\} }}+2\,\Delta_{{ \left\{ 3,4,5,7
 \right\} }}{q}^{2}-
 \]
 \[
 -\Delta_{{ \left\{ 2,4,7,8 \right\} }}{q}^{2}+
\Delta_{{ \left\{ 2,3,7,8 \right\} }}{q}^{2}-\Delta_{{ \left\{ 1,4,7,8
 \right\} }}p-3\,\Delta_{{ \left\{ 3,4,5,7 \right\} }}q)
 \]
Assuming $\left\{ \Delta_{{ \left\{ 1,4,7,8 \right\} }},\Delta_{{ \left\{ 2,3,7
,8 \right\} }},\Delta_{{ \left\{ 2,4,7,8 \right\} }},\Delta_{{
 \left\{ 3,4,5,7 \right\} }},\Delta_{{ \left\{ 3,4,7,8 \right\} }}
 \right\}=\{1\}$, we can try to resolve all remaining (Plu"cker) relations.
In this case
\begin{equation}
{\it Fset}=\left[-{\frac {2\,q-qp}{qp}},1,1,1,1,1,-{\frac {2\,qp-2\,q+p}{p
 \left( -1+2\,q \right) }},{\frac {-2\,q+3\,qp-2\,{q}^{2}p+2\,{q}^{2}}
{pq \left( -1+2\,q \right) }},1\right]
\end{equation} 
Now we're almost ready to write our target polynomial. The only thing left is to convert
the indices ${i_1,i_2,\dots,i_s}$ to the corresponding Young diagram. In order to do this,
expand the following determinant by Cauchy-Binet formula:
We have:
\begin{equation}
\Delta_{\{1,3,7,8\}}\simeq\left| \begin {array}{cccc} h_{{2}}&h_{{4}}&h_{{5}}&h_{{6}}
\\\noalign{\medskip}h_{{1}}&h_{{3}}&h_{{4}}&h_{{5}}
\\\noalign{\medskip}0&h_{{2}}&h_{{3}}&h_{{4}}\\\noalign{\medskip}0&h_{
{1}}&h_{{2}}&h_{{3}}\end {array} \right|=s_{\tiny\yng(2,3,3,3)},  
\Delta_{\{1,4,7,8\}}\simeq 
\left| \begin {array}{cccc} h_{{2}}&h_{{3}}&h_{{5}}&h_{{6}}
\\\noalign{\medskip}h_{{1}}&h_{{2}}&h_{{4}}&h_{{5}}
\\\noalign{\medskip}0&h_{{1}}&h_{{3}}&h_{{4}}\\\noalign{\medskip}0&0&h
_{{2}}&h_{{3}}\end {array} \right|=s_{\tiny\yng(2,2,3,3)}, \Delta_{\{2,3,7,8\}}\simeq s_{\tiny\yng(1,3,3,3)},
\end{equation}
\[
\Delta_{\{2,4,7,8\}}\simeq s_{\tiny\yng(1,2,3,3)}, \Delta_{\{3,4,5,7\}}\simeq s_{\tiny\yng(1,1,4,5)},
\Delta_{\{3,4,5,8\}}\simeq s_{\tiny\yng(1,1,4,4)},\Delta_{\{3,4,6,7\}}\simeq s_{\tiny\yng(1,1,3,5)},
\Delta_{\{3,4,6,8\}}\simeq s_{\tiny\yng(1,1,3,4)}, \Delta_{\{3,4,7,8\}}\simeq s_{\tiny\yng(1,1,3,3)}.
\]

\subsection{$Q=A_2,d_i=6,e_i=3, \mu=\sigma(p,q,6)$}
Take 2 vertices $v_1,v_2$ with $V_1=V_2=\mathbb{C}^6$, 
\begin{equation}
I=\left\{\underbrace{ 4,5,6},\underbrace{10,11,12} \right\} 
\end{equation}
Recall that each $\Delta$ is a $(6,6)$ complex minor.
The following minors are supposed to vanish:
\begin{equation}
\{ \Delta_{{ \left\{ 1,4,5,6,10,11 \right\} }},\Delta_{{
 \left\{ 1,4,5,6,10,12 \right\} }},\Delta_{{ \left\{ 1,4,5,6,11,12
 \right\} }},\Delta_{{ \left\{ 2,4,5,6,10,11 \right\} }},\Delta_{{
 \left\{ 2,4,5,6,10,12 \right\} }},
\end{equation}
\[
\Delta_{{ \left\{ 2,4,5,6,11,12
 \right\} }},\Delta_{{ \left\{ 3,4,5,6,10,11 \right\} }},\Delta_{{
 \left\{ 3,4,5,6,10,12 \right\} }},\Delta_{{ \left\{ 3,4,5,6,11,12
 \right\} }} \} =\{0\}
\]
Total $\#QR=3^3$ . (3 for each $i in \{4,5,6\}$: $F(v,i,j)$).
For arbitrary $\mu$ these relations look really huge!

The choice $\mu=Id(6)$ kills most of the $Fset$ entries, so we can write 
$Fset$ as $(12,12)$ sparse matrix with $\{0,\pm1\}$ entries of the shape
(empty space is filled with zeroes):
\begin{verbatim}
Fset=Matrix([seq(seq(Fset[1+i .. 12+i], i = 12*j), j = 0 .. 11)])
\end{verbatim}
\begin{equation}
Fset={\tiny{ \left( \begin {array}{cccccccccccc} &&&1&-1&1&&&&&&
\\\noalign{\medskip}&&&1&-1&1&&&&&&\\\noalign{\medskip}&&&
1&-1&1&&&&&&\\\noalign{\medskip}&&&1&&&&&&&&
\\\noalign{\medskip}&&&&1&&&&&&&\\\noalign{\medskip}&&&
&&1&&&&&&\\\noalign{\medskip}&&&&&&&&&1&1&1
\\\noalign{\medskip}&&&&&&&&&1&1&1\\\noalign{\medskip}&&&
&&&&&&1&1&1\\\noalign{\medskip}&&&&&&&&&1&&
\\\noalign{\medskip}&&&&&&&&&&1&\\\noalign{\medskip}&&&
&&&&&&&&1\end {array} \right)}}
\end{equation}
Finally, convert all survived $Fset$ entries to corresponding Young diagrams to get out
target polynomial:
\begin{equation}
\tilde{P}(A_2,6,Id)=s_{\tiny\yng(2,2,2,4,4,4)}-s_{\tiny\yng(2,2,3,4,4,4)}+s_{\tiny\yng(2,2,4,4,4,4)}+
s_{\tiny\yng(1,2,2,4,4,4)}-s_{\tiny\yng(1,2,3,4,4,4)}+s_{\tiny\yng(1,2,4,4,4,4)}+\dots,
\end{equation}
total $\#(0,+,-)=(120,21,3)$

\subsection{$Q=A_3,d_i=2,\mu_{i=1..2}=Mat(2,2,\mathbb{C}),\overline{e}=(1,1,1)\}$}
Now let's add one more vertex:
\begin{equation}
Q=A_3:\quad \circ \xleftarrow{\mu_{1}} \circ \xleftarrow{\mu_2} \circ, \quad \mu_i
\end{equation}
\begin{equation}
Gr_{(1,1,1)}(A_3,M)=
\end{equation}
Choose $I=\{2,4,6\}$, so now we have the embedding
\begin{equation}
Gr_{(1,1,1)}(A_3,M)\longrightarrow Gr(3,6) \longrightarrow \mathbb{P}^{6+3-1}
\end{equation}
For $v=1$ we have the following non-trivial quiver relations:
\begin{equation}
F(1,2,3)=\mu_{{1,1,4}}\Delta_{{ \left\{ 1,4,6 \right\} }}\Delta_{{ \left\{ 2,3,
6 \right\} }}+\mu_{{1,2,4}}\Delta_{{ \left\{ 2,4,6 \right\} }}\Delta_{
{ \left\{ 2,3,6 \right\} }}-\mu_{{1,1,3}}\Delta_{{ \left\{ 1,4,6
 \right\} }}\Delta_{{ \left\{ 2,4,6 \right\} }}-\mu_{{1,2,3}}{\Delta_{
{ \left\{ 2,4,6 \right\} }}}^{2}=0
\end{equation}
Note that if $\mu=Id$, the solution looks as follows:
\begin{equation}
Fset= \left( \begin {array}{cccccc} &\Delta_{{ \left\{ 1,4,6 \right\} }}&
&\Delta_{{ \left\{ 1,2,6 \right\} }}&&\Delta_{{ \left\{ 1,2,4
 \right\} }}\\\noalign{\medskip}&1&&&&\\\noalign{\medskip}&&&
\Delta_{{ \left\{ 1,4,6 \right\} }}&&\\\noalign{\medskip}&&&1&&
\\\noalign{\medskip}&&&&&\Delta_{{ \left\{ 1,4,6 \right\} }}
\\\noalign{\medskip}&&&&&1\end {array} \right),
\end{equation}
where all empty cells are zeroes.

\subsection{$Q=A_3, d_i=4, \mu_1=\sigma(p,q,4),\mu_2=Id$}
Initial data:
\begin{equation}
\mathcal{B}=\{1,2,\dots,12\}, I=\{3,4,7,8,11,12\}
\end{equation}
\begin{equation}
 \mu_{{1}}= \left(\begin {array}{cccc} p&1-p&0&0
\\\noalign{\medskip}1-p&p&0&0\\\noalign{\medskip}0&0&q&1-q
\\\noalign{\medskip}0&0&1-q&q\end {array} \right) ,\mu_{{2}}= \left( 
\begin {array}{cccc} 1&0&0&0\\\noalign{\medskip}0&1&0&0
\\\noalign{\medskip}0&0&1&0\\\noalign{\medskip}0&0&0&1\end {array}
 \right)
\end{equation}
This case reads
\begin{equation}
{\it Fset'}=\left\{1,\Delta_{{ \left\{ 2,3,4,7,8,11 \right\} }},\Delta_{{
 \left\{ 2,3,4,7,8,12 \right\} }},\Delta_{{ \left\{ 3,4,7,8,9,11
 \right\} }},\Delta_{{ \left\{ 3,4,7,8,10,11 \right\} }},-{\frac {-q+2
\,qp}{q \left( -1+2\,p \right) }},-{\frac {2\,{p}^{2}q-qp}{qp \left( -
1+2\,p \right) }}\right\}
\end{equation}

%

\section{Diagrammatic notation and exact formulas for invariants}

\textbf{$\tilde{P}(A_n,2k,\sigma)$:}
\begin{equation}
\tilde{P}(A_2,4,\sigma)={\frac { \left( -2\,q-qp \right)  \left( {h_{{3}}}^{3}h_{{2}}-{h_{{3}}
}^{2}h_{{1}}h_{{4}}-2\,h_{{3}}h_{{4}}{h_{{2}}}^{2}+2\,{h_{{4}}}^{2}h_{
{2}}h_{{1}}-h_{{4}}h_{{5}}{h_{{1}}}^{2}+h_{{5}}{h_{{2}}}^{3}-h_{{6}}{h
_{{2}}}^{2}h_{{1}}+h_{{6}}h_{{3}}{h_{{1}}}^{2} \right) }{qp}}+
\end{equation}
\[
+{h_{{3}}
}^{2}{h_{{2}}}^{2}-h_{{5}}h_{{3}}{h_{{1}}}^{2}-h_{{4}}{h_{{2}}}^{3}+2
\,h_{{5}}{h_{{2}}}^{2}h_{{1}} 
-h_{{6}}h_{{2}}{h_{{1}}}^{2}-2\,h_{{3}}h_{{4}}h_{{2}}h_{{1}}+2\,{h_{{4
}}}^{2}{h_{{1}}}^{2}+{h_{{3}}}^{2}h_{{2}}h_{{1}}-h_{{3}}h_{{4}}{h_{{1}
}}^{2}+
\]
\[
{\frac { \left( -2\,qp-2\,q+p \right)  \left( h_{{5}}h_{{3}}{h_{{1}}}^
{2}-h_{{6}}h_{{2}}{h_{{1}}}^{2} \right) }{p \left( 2\,q-1 \right) }}+{
\frac { \left( -2\,q+3\,qp-2\,{q}^{2}p+2\,{q}^{2} \right)  \left( h_{{
3}}h_{{4}}{h_{{1}}}^{2}-h_{{5}}h_{{2}}{h_{{1}}}^{2} \right) }{pq
 \left( 2\,q-1 \right) }}-
\]
\[
-h_{{4}}{h_{{2}}}^{2}h_{{1}}+h_{{5}}h_{{2}}{h_{{1}}}^{2}+h_{{4}}
h_{{5}}{h_{{1}}}^{2}-h_{{6}}h_{{3}}{h_{{1}}}^{2}+{h_{{3}}}^{2}{h_{{1}}}^{2}-h_{{4}}h_{{2}}{h_{{1}}}^{2}
\]

What is the asymptotic formula for $\#+,-$ in $Fset(Id,m)$, when $m\to \infty$ ? The first few terms are
shown on the figure:
\begin{figure}[h!]
\centering
\includegraphics[scale=0.3]{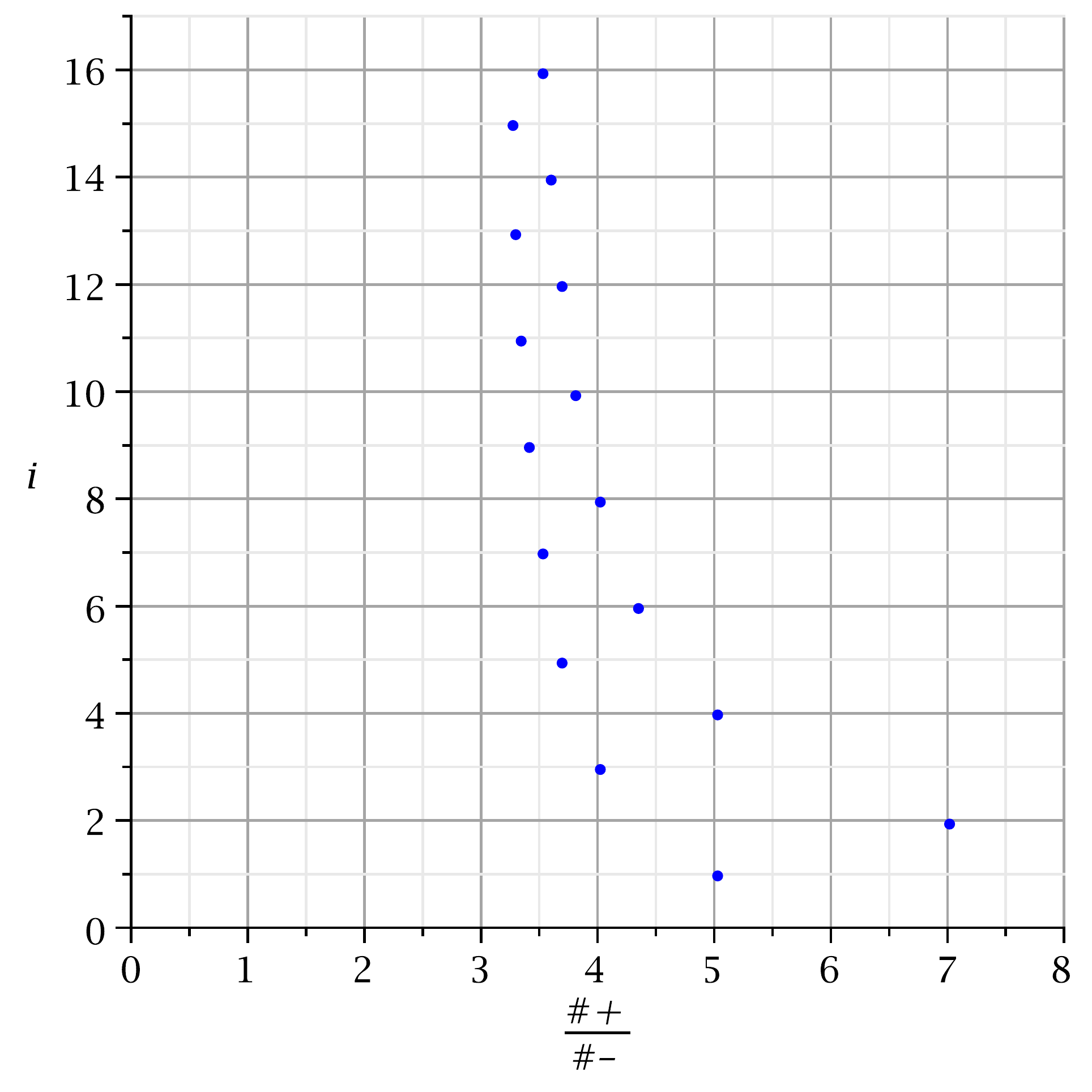}
\caption{Convergence of $\{\frac{\#+}{\#-}\}$ sequence, $A_2$ quiver, $\mu=Id,\mathrm{dim}(M)=i$}
\end{figure}
%
%
$Fset$ matrix arrangement for $A_4$ quiver with $d_i=4,\mu_i=Id$:
\begin{equation}
Fset={ \left( \begin {array}{cccccccccccccccc} && x_{1}& x_{2}&&&
x_{3}& x_{4}&&& x_{5}& x_{6}&&& x_{7}& x_{8}
\\\noalign{\medskip}&& x_{1}& x_{11}&&& x_{12}& x_{13}&&
& x_{14}& x_{15}&&& x_{16}& x_{17}\\\noalign{\medskip}&&1&
&&&&&&&&&&&&\\\noalign{\medskip}&&&1&&&&&&&&&&
&&\\\noalign{\medskip}&&&&&& x_{1}&- x_{2}&&& x_{18}&
 x_{19}&&& x_{2}& x_{21}\\\noalign{\medskip}&&&&&& 
x_{1}&- x_{11}&&& x_{22}& x_{23}&&& x_{24}& x_{25}
\\\noalign{\medskip}&&&&&&1&&&&&&&&&
\\\noalign{\medskip}&&&&&&&1&&&&&&&&
\\\noalign{\medskip}&&&&&&&&&& x_{1}& x_{2}&&&&
\\\noalign{\medskip}&&&&&&&&&& x_{1}& x_{11}&&&&
\\\noalign{\medskip}&&&&&&&&&&1&&&&&
\\\noalign{\medskip}&&&&&&&&&&&1&&&&
\\\noalign{\medskip}&&&&&&&&&&&&&& x_{1}&- x_{2}
\\\noalign{\medskip}&&&&&&&&&&&&&& x_{1}&- x_{11}
\\\noalign{\medskip}&&&&&&&&&&&&&&1&
\\\noalign{\medskip}&&&&&&&&&&&&&&&1\end {array}
 \right)}
\end{equation}
where $x_j=\Delta_{I(j)}$. Compare this to $A_2,d_i=8,\mu_i=Id$ case:
\begin{equation}
\tiny
{\left( \begin {array}{cccccccccccccccc} &&&&x_{{1}}&-x_{{2}}&x_{{
3}}&-x_{{4}}&&&&&&&&\\\noalign{\medskip}&&&&x_{{5}}&x_{{6}
}&x_{{7}}&x_{{8}}&&&&&&&&\\\noalign{\medskip}&&&&x_{{9}}&x
_{{1}}&x_{{11}}&x_{{12}}&&&&&&&&\\\noalign{\medskip}&&&&x
_{{13}}&x_{{14}}&x_{{15}}&x_{{16}}&&&&&&&&\\\noalign{\medskip}
&&&&1&&&&&&&&&&&\\\noalign{\medskip}&&&&&1&&&&&
&&&&&\\\noalign{\medskip}&&&&&&1&&&&&&&&&
\\\noalign{\medskip}&&&&&&&1&&&&&&&&
\\\noalign{\medskip}&&&&&&&&&&&&x_{{1}}&x_{{2}}&x_{{3}}&x_
{{4}}\\\noalign{\medskip}&&&&&&&&&&&&x_{{5}}&-x_{{6}}&x_{{
7}}&-x_{{8}}\\\noalign{\medskip}&&&&&&&&&&&&x_{{9}}&-x_{{
1}}&x_{{11}}&-x_{{12}}\\\noalign{\medskip}&&&&&&&&&&&&x_{
{13}}&-x_{{14}}&x_{{15}}&-x_{{16}}\\\noalign{\medskip}&&&&&&&&
&&&&1&&&\\\noalign{\medskip}&&&&&&&&&&&&&1&&
\\\noalign{\medskip}&&&&&&&&&&&&&&1&
\\\noalign{\medskip}&&&&&&&&&&&&&&&1\end {array}
 \right)}
\end{equation}

\underline{Totally positive} $Fset$ arrangement for $A_4 \times {\mathrm{mutations}},d_i=6,\mu=Id$:
\begin{equation}
 \left( \begin {array}{cccccccccccccccccccccccc} &&&1&
&1&&&&1&1&1&&&&1&1&1&&&&1&1&1
\\\noalign{\medskip}&&&1&&1&&&&1&1&1&&
&&1&1&1&&&&1&1&1\\\noalign{\medskip}&&
&1&&1&&&&1&1&1&&&&1&1&1&&&&1&1
&1\\\noalign{\medskip}&&&1&&&&&&&
&&&&&&&&&&&&&
\\\noalign{\medskip}&&&&1&&&&&&
&&&&&&&&&&&&&
\\\noalign{\medskip}&&&&&1&&&&&
&&&&&&&&&&&&&
\\\noalign{\medskip}&&&&&&&&&1&
&1&&&&1&1&1&&&&1&1&1\\\noalign{\medskip}
&&&&&&&&&1&&1&&&&1
&1&1&&&&1&1&1\\\noalign{\medskip}&&&&&
&&&&1&&1&&&&1&1&1&&&&1&1&1
\\\noalign{\medskip}&&&&&&&&&1&
&&&&&&&&&&&&&
\\\noalign{\medskip}&&&&&&&&&&
1&&&&&&&&&&&&&
\\\noalign{\medskip}&&&&&&&&&&
&1&&&&&&&&&&&&
\\\noalign{\medskip}&&&&&&&&&&
&&&&&1&&1&&&&1&1&1
\\\noalign{\medskip}&&&&&&&&&&
&&&&&1&&1&&&&1&1&1
\\\noalign{\medskip}&&&&&&&&&&
&&&&&1&&1&&&&1&1&1
\\\noalign{\medskip}&&&&&&&&&&
&&&&&1&&&&&&&&
\\\noalign{\medskip}&&&&&&&&&&
&&&&&&1&&&&&&&
\\\noalign{\medskip}&&&&&&&&&&
&&&&&&&1&&&&&&
\\\noalign{\medskip}&&&&&&&&&&
&&&&&&&&&&&1&&1
\\\noalign{\medskip}&&&&&&&&&&
&&&&&&&&&&&1&&1
\\\noalign{\medskip}&&&&&&&&&&
&&&&&&&&&&&1&&1
\\\noalign{\medskip}&&&&&&&&&&
&&&&&&&&&&&1&&
\\\noalign{\medskip}&&&&&&&&&&
&&&&&&&&&&&&1&
\\\noalign{\medskip}&&&&&&&&&&
&&&&&&&&&&&&&1
\end {array} \right)
\end{equation}
\section{•}
\scalebox{0.5}{\parbox{\linewidth}{
\begin{equation}
\begin{aligned}
(2,2)=& \left( s_{\tiny\yng(2,2)}+s_{\tiny\yng(1,3)} \right) \Delta_{\tiny\yng(2,2)} & jsh=\{{\tiny\yng(1)}\} \\
(3,4)=&\left( s_{\tiny\yng(2,2,3,3,5,5)}+s_{\tiny\yng(1,1,4,4,5,5)}+s_{\tiny\yng(1,1,3,3,6,6)} \right) 
\Delta_{\tiny\yng(1,1,3,3,6,6)}+ \left( s_{\tiny\yng(1,2,3,3,5,5)}+s_{\tiny\yng(1,1,3,4,5,5)}+s_{
\tiny\yng(1,1,3,3,5,6)} \right) \Delta_{\tiny\yng(1,2,3,3,5,5)} & jsh=\{{\tiny\yng(1),\yng(1,1)}\} \\
(4,6)=& 
 \left( s_{\tiny\yng(1,1,2,4,4,4,7,7,7,10,10,10)}+s_{\tiny\yng(1,1,1,4,4,5,7,7,7,10,10,10)}
+s_{\tiny\yng(1,1,1,4,4,4,7,7,8,10,10,10)}+s_{\tiny\yng(1,1,1,4,4,4,7,7,7,10,10,11)}
 \right) \Delta_{\tiny\yng(1,1,1,4,4,5,7,7,7,10,10,10)}+ \left( s_{\tiny\yng(1,1,1,4,4,7,7,7,7,10,10,10)}+
s_{\tiny\yng(1,1,1,4,4,4,7,7,7,10,10,13)}+s_{\tiny\yng(1,1,1,4,4,4,7,7,10,10,10,10)}+
s_{\tiny\yng(1,1,4,4,4,4,7,7,7,10,10,10)} \right) \Delta_{\tiny\yng(1,1,1,4,4,7,7,7,7,10,10,10)}+ \\
&
\left( s_{\tiny\yng(1,2,2,4,4,4,7,7,7,10,10,10)}+s_{\tiny\yng(1,1,1,4,4,4,7,7,7,10,11,11)}+
s_{\tiny\yng(1,1,1,4,5,5,7,7,7,10,10,10)}+s_{\tiny\yng(1,1,1,4,4,4,7,8,8,10,10,10)} \right) 
\Delta_{\tiny\yng(1,1,1,4,5,5,7,7,7,10,10,10)}+
 \left( s_{\tiny\yng(1,1,1,4,4,4,7,8,10,10,10,10)}+s_{\tiny\yng(1,2,4,4,4,4,7,7,7,10,10,10)}+
s_{\tiny\yng(1,1,1,4,5,7,7,7,7,10,10,10)}+s_{\tiny\yng(1,1,1,4,4,4,7,7,7,10,11,13)} \right)
 \Delta_{\tiny\yng(1,1,1,4,5,7,7,7,7,10,10,10)}+ \\
 &
 \left( s_{\tiny\yng(2,2,2,4,4,4,7,7,7,10,10,10)}+
s_{\tiny\yng(1,1,1,4,4,4,8,8,8,10,10,10)}+s_{\tiny\yng(1,1,1,5,5,5,7,7,7,10,10,10)}+
s_{\tiny\yng(1,1,1,4,4,4,7,7,7,11,11,11)} \right) \Delta_{\tiny\yng(1,1,1,5,5,5,7,7,7,10,10,10)}+ 
\left( s_{\tiny\yng(2,2,4,4,4,4,7,7,7,10,10,10)}+s_{\tiny\yng(1,1,1,5,5,7,7,7,7,10,10,10)}+
s_{\tiny\yng(1,1,1,4,4,4,7,7,7,11,11,13)}+s_{\tiny\yng(1,1,1,4,4,4,8,8,10,10,10,10)} \right)
 \Delta_{\tiny\yng(1,1,1,5,5,7,7,7,7,10,10,10)}
& jsh=\{{\tiny\yng(1),\yng(1,1),\yng(1,3),\yng(3),\yng(1,1,1),\yng(1,1,3)}\} \\
(5,8)=&
\left( s_{\tiny\yng(1,1,1,1,5,5,5,5,9,9,9,10,13,13,13,13,17,17,17,17)}+
s_{\tiny\yng(1,1,1,2,5,5,5,5,9,9,9,9,13,13,13,13,17,17,17,17)}+
s_{\tiny\yng(1,1,1,1,5,5,5,6,9,9,9,9,13,13,13,13,17,17,17,17)}+
s_{\tiny\yng(1,1,1,1,5,5,5,5,9,9,9,9,13,13,13,14,17,17,17,17)}+
s_{\tiny\yng(1,1,1,1,5,5,5,5,9,9,9,9,13,13,13,13,17,17,17,18)} \right) 
\Delta_{\tiny\yng(1,1,1,1,5,5,5,5,9,9,9,9,13,13,13,14,17,17,17,17)}+\dots 
& jsh=\{{\tiny\yng(1),\yng(1,1),\dots,\yng(1,1,1,1),\yng(1,1,1,3)}\}
\\
\end{aligned}
\end{equation}
}}

\scalebox{.2}{
\tiny\yng(1,1,1,2,2,5,7,7,7,7,7,7,13,13,13,13,13,13,19,19,19,19,19,19,25,25,25,25,25,25,31,31,31,31,31,31,37,37,37,37,37,37,43,43,43,43,43,43,49,49,49,49,49,49,55,55,55,55,55,55,61,61,61,61,61,61,67,67,67,67,67,67)}

$\frac{1}{2}(\tilde{P}_{+}-\tilde{P}_{-})$ for cyclic 6-vertex graph:

\scalebox{0.5}{\parbox{\linewidth}{
\begin{equation}
 \left( -s_{\tiny\yng(1,1,1,4,4,4,7,7,7,10,10,10,13,13,16,16,16,16)}+
s_{\tiny\yng(1,1,1,4,4,4,7,7,7,10,10,13,13,13,13,16,16,16)} \right)
 \Delta_{\tiny\yng(1,1,1,4,4,4,7,7,7,10,10,13,13,13,13,16,16,16)}+ 
\left( s_{\tiny\yng(1,1,1,4,4,4,7,7,7,10,11,13,13,13,13,16,16,16)}-
s_{\tiny\yng(1,1,1,4,4,4,7,7,7,10,10,10,13,14,16,16,16,16)} \right)
 \Delta_{\tiny\yng(1,1,1,4,4,4,7,7,7,10,11,13,13,13,13,16,16,16)}+
 \left( -s_{\tiny\yng(1,1,1,4,4,4,7,7,7,10,10,10,14,14,16,16,16,16)}+
s_{\tiny\yng(1,1,1,4,4,4,7,7,7,11,11,13,13,13,13,16,16,16)} \right) 
\Delta_{\tiny\yng(1,1,1,4,4,4,7,7,7,11,11,13,13,13,13,16,16,16)}
\end{equation}
}}
\newpage
\textbf{Lemma 1.} \emph{Asymptotic formula for $\tilde{P}_{\mathrm{max}}(Q_{\mathrm{univ}}(n),n,Id)$:}
\begin{itemize}
\item \emph{$m \times m$ ``cubic ladder'' decomposition $\simeq \Delta$-summands in $\tilde{P_{I}}$}
\item \emph{each $\Delta$-summand consists of Young diagrams factorized by the ``ladder-block'' addition}
\end{itemize}

The formula:
\begin{equation}
\boxed{\tilde{P}_{\mathrm{max}}(Q(m),n,Id)=\sum\limits_{i=1..m\cdot\mathrm{cl}\left(\frac{m}{2}\right)}\sum\limits_{j=1..n}s_{\lambda(i,j)}\Delta_{\lambda(i,*)},}
\end{equation}

$\frac{1}{2}(\tilde{P}_{+}-\tilde{P}_{-})$ for cyclic 8-vertex graph:

\scalebox{0.22}{\parbox{\linewidth}{
\begin{equation}
\begin{aligned}
\left( -s_{\tiny\yng(1,1,1,1,5,5,5,5,9,9,9,9,13,13,13,13,18,18,18,20,21,21,21,21,25,25,25,25,29,29,29,29)}+s_{\tiny\yng(1,1,1,1,5,5,5,5,9,9,9,9,13,13,13,13,17,17,17,17,22,22,22,24,25,25,25,25,29,29,29,29)} \right) \Delta_{\tiny\yng(1,1,1,1,5,5,5,5,9,9,9,9,13,13,13,15,17,17,17,17,21,21,21,21,25,25,25,25,29,29,29,29)}+ \left( s_{\tiny\yng(1,1,1,1,5,5,5,5,9,9,9,9,13,13,13,13,17,17,17,17,21,21,21,21,26,26,26,28,29,29,29,29)}-s_{\tiny\yng(1,1,1,1,5,5,5,5,9,9,9,9,13,13,13,13,18,18,18,18,21,21,21,21,25,25,25,25,29,29,29,29)} \right) 
\Delta_{\tiny\yng(1,1,1,1,5,5,5,5,9,9,9,9,13,13,13,17,17,17,17,17,21,21,21,21,25,25,25,25,29,29,29,29)}+
 \left( s_{\tiny\yng(1,1,1,1,5,5,5,5,9,9,9,9,13,13,13,13,17,17,17,17,21,21,21,24,25,25,25,25,29,29,29,29)}-s_{\tiny\yng(1,1,1,1,5,5,5,5,9,9,9,9,13,13,13,13,17,17,17,20,21,21,21,21,25,25,25,25,29,29,29,29)} \right) \Delta_{\tiny\yng(1,1,1,1,5,5,5,5,9,9,9,9,13,13,14,15,17,17,17,17,21,21,21,21,25,25,25,25,29,29,29,29)}+
  \left( s_{\tiny\yng(1,1,1,1,5,5,5,5,9,9,9,9,13,13,13,13,17,17,17,17,21,21,21,21,25,25,25,28,29,29,29,29)}-
 s_{\tiny\yng(
1,1,1,1,5,5,5,5,9,9,9,9,13,13,13,13,17,17,17,18,21,21,21,21,25,25,25,25,29,29,29,29)} \right) \Delta_{\tiny\yng(1,1,1,1,5,5,5,5,9,9,9,9,13,13,14,17,17,17,17,17,21,21,21,21,25,25,25,25,29,29,29,29)}+
&
\\
 \left( s_{\tiny\yng(1,1,1,1,5,5,5,5,9,9,9,9,13,13,13,13,17,17,17,17,21,21,22,24,25,25,25,25,29,29,29,29)}-s_{\tiny\yng(1,1,1,1,5,5,5,5,9,9,9,9,13,13,13,13,17,17,18,20,21,21,21,21,25,25,25,25,29,29,29,29)} \right) \Delta_{\tiny\yng(1,1,1,1,5,5,5,5,9,9,9,9,13,14,14,15,17,17,17,17,21,21,21,21,25,25,25,25,29,29,29,29)}+ \left( 
s_{\tiny\yng(1,1,1,1,5,5,5,5,9,9,9,9,13,13,13,13,17,17,17,17,21,21,21,21,25,25,26,28,29,29,29,29)}-s_{\tiny\yng(1,1,1,1,5,5,5,5,9,9,9,9,13,13,13,13,17,17,18,18,21,21,21,21,25,25,25,25,29,29,29,29)} \right) \Delta_{\tiny\yng(1,1,1,1,5,5,5,5,9,9,9,9,13,14,14,17,17,17,17,17,21,21,21,21,25,25,25,25,29,29,29,29)}+
 \left( s_{\tiny\yng(1,1,1,1,5,5,5,5,9,9,9,9,13,13,13,13,17,17,17,17,21,22,22,24,25,25,25,25,29,29,29,29)}-s_{\tiny\yng(1,1,1,1,5,5,5,5,9,9,9,9,13,13,13,13,17,18,18,20,21,21,21,21,25,25,25,25,29,29,29,29)} \right) \Delta_{\tiny\yng(1,1,1,1,5,5,5,5,9,9,9,9,14,14,14,15,17,17,17,17,21,21,21,21,25,25,25,25,29,29,29,29)}+
 \left( s_{\tiny\yng(1,1,1,1,5,5,5,5,9,9,9,9,13,13,13,13,17,17,17,17,21,21,21,21,25,26,26,28,29,29,29,29)}-s_{\tiny\yng(1,1,1,1,5,5,5,5,9,9,9,9,13,13,13,13,17,18,18,18,21,21,21,21,25,25,25,25,29,29,29,29)}
 \right) \Delta_{\tiny\yng(1,1,1,1,5,5,5,5,9,9,9,9,14,14,14,17,17,17,17,17,21,21,21,21,25,25,25,25,29,29,29,29)}
\end{aligned}
\end{equation}
}}
\newpage
\section{Re-evaluation of $\tilde{P}$ for Coxeter quivers. ``Vertex at $\infty$''}
Let $Q$ be a quiver with $m$ vertices and all possible arrows (but only 2 opposite arrows for each pair of vertices
is allowed, so the direction $v=\pm 1$) .
Choose $n,e\in \mathbb{N}$ and the integer-valued vectors
\begin{equation}
B=\{1,\dots,m\cdot{}n\}, I=\{jn-i+1\}_{i=1..[n/2],j=1..m}\subset B,
|I|=e,
\end{equation}
where $[x] \text{ is the integer part of } x$.
Fix the representation $R(Q)$: vertices $\sim V_i(\mathbb{Z}^n)$, arrows $\sim$ transposition matrices $M_{v,i,j}=\sigma(\lambda_{i,j})^{\mathrm{sign}(v)}$. 

Each $\sigma(\lambda_{i,j})$ is a $(n \times n)$ transposition matrix, partitioned by $\lambda_{i,j}$ for the $(i,j)$-th arrow.
We fix \emph{the particular choice} $\boxed{R=\Sigma_0}$:
\begin{equation}
(i,j)\rightarrow\lambda_{i,j}: (1,4)\rightarrow {\tiny\yng(1,2,3,4)}, (1,3)\rightarrow{\tiny\yng(1,2,3,3)},  (1,2)\rightarrow{\tiny\yng(1,2,2,2)},  (2,4)\rightarrow{\tiny\yng(2,2,3,4)}, (2,3)\rightarrow{\tiny\yng(2,2,3,3)},\dots
\end{equation}
Assume also $\mathrm{dim(Gr}(Q,R))=\boxed{\overline{e}=\left[\frac{1}{2}\overline{n}\right]}$ for the corresponding quiver grassmannian \cite{2,4}, $e_i=e_j\ \forall i,j$.
For instance, the case $(m,n)=(4,4)$ is given by
\begin{equation}
\mathbf{M}_{+}={\tiny  \left( \begin {array}{cccc} 0& \left[ \begin {array}{cccc} 1&0&0&0
\\\noalign{\medskip}0&0&1&0\\\noalign{\medskip}0&0&0&1
\\\noalign{\medskip}0&1&0&0\end {array} \right] & \left[ 
\begin {array}{cccc} 1&0&0&0\\\noalign{\medskip}0&1&0&0
\\\noalign{\medskip}0&0&0&1\\\noalign{\medskip}0&0&1&0\end {array}
 \right] & \left[ \begin {array}{cccc} 1&0&0&0\\\noalign{\medskip}0&1&0
&0\\\noalign{\medskip}0&0&1&0\\\noalign{\medskip}0&0&0&1\end {array}
 \right] \\\noalign{\medskip}0&0& \left[ \begin {array}{cccc} 0&1&0&0
\\\noalign{\medskip}1&0&0&0\\\noalign{\medskip}0&0&0&1
\\\noalign{\medskip}0&0&1&0\end {array} \right] & \left[ 
\begin {array}{cccc} 0&1&0&0\\\noalign{\medskip}1&0&0&0
\\\noalign{\medskip}0&0&1&0\\\noalign{\medskip}0&0&0&1\end {array}
 \right] \\\noalign{\medskip}0&0&0& \left[ \begin {array}{cccc} 0&1&0&0
\\\noalign{\medskip}0&0&1&0\\\noalign{\medskip}1&0&0&0
\\\noalign{\medskip}0&0&0&1\end {array} \right] \\\noalign{\medskip}0&0
&0&0\end {array} \right)
 }
\end{equation}
and $\mathbf{M}_{-}$ equals by-block inverse of $\mathbf{M}_{+}$.
Choose some path $x\in \mathcal{P}(Q)$ and solve the quiver relations with vanishing conditions \cite{4}:
\begin{equation}
\begin{aligned}
\quad &\sum_{s',t'} (-1)_{\epsilon(s',t')}M_{v,s',t'} \ \Delta_s^{s'} \ \Delta_{t'}^t 
- \ \sum_{s'} (-1)_{\epsilon(s',t)} M_{v,s',t} \ \Delta_{s}^{s'} \ \Delta_I = 0, \\
& \Delta_J=0 \text{ for some non-proper } J\text{'s}
\end{aligned}
\end{equation}
Define symmetric polynomials
\begin{equation}
\boxed{\tilde{P}_R(x)=\sum\limits_{i,j}\Delta_{\mu_{i,i}}s_{\mu_{i,j}}(p_1,p_2,\dots)},
\end{equation}
where $\mu_{i,j}$ is a partition derived from the Cauchy-Binet expansion for each $\Delta$'s, $s_{\mu}$ is a Schur function.

\textbf{Prop 1.}
\emph{The map $f:x \mapsto \tilde{P}_R(x)$ is a homomorphism on a ``proper'' subset $S\subseteq\mathcal{P}$, s.t.
for any $x,y \in S$ }
\begin{equation}
\tilde{P}_R(x\cdot y)=\frac{1}{2}\sum_{\mu,\eta,\nu}\left(\mathrm{coeff}(\tilde{P}_R(x),\Delta_{\mu})-\mathrm{coeff}(\tilde{P}_R(y),\Delta_{\eta})\right)\Delta_{\nu}\simeq_{|_{Res(\Delta)}}\frac{1}{2}\left( \tilde{P}_R(x)-\tilde{P}_R(y) \right)
\end{equation}
\underline{Note:} this construction is very sensitive to the choice of $R$ (if we swap the transposition matrices, $S$ will change).
%
%
\begin{figure}[h!]
\centering
\includegraphics[scale=0.35]{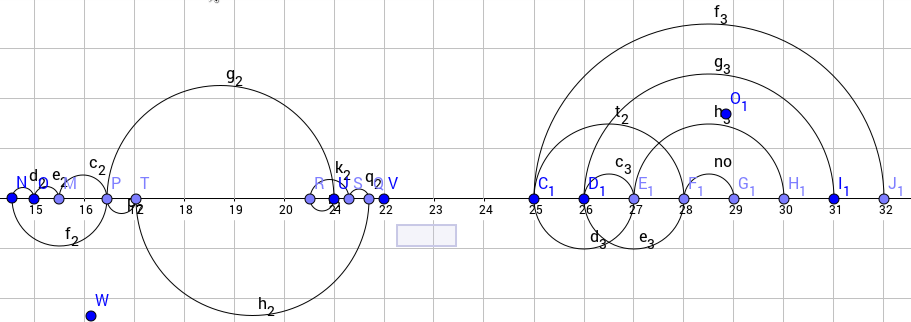}
\caption{Examples of ``good'' (left) and ``bad'' (right) combinations of paths in $S\subseteq\mathcal{P}$
for $(m,n)=(10,10)$.
``Good'' means the image of $f$ is not empty}
\end{figure}
\begin{figure}[h!]
\centering
\includegraphics[scale=0.3]{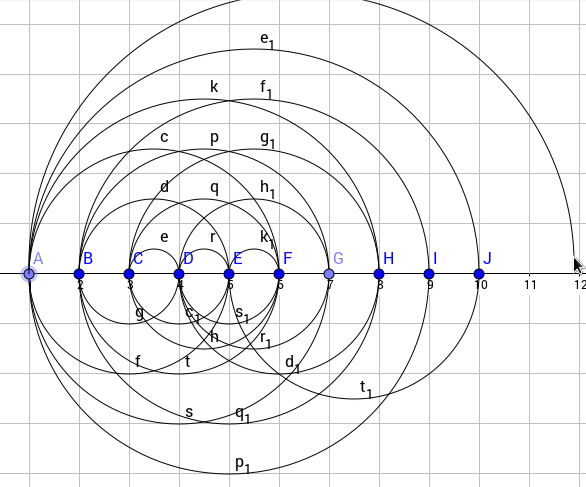}
\caption{Examples of a non-trivial loops \emph{spiral type 1} (Coxeter odd type $\tilde{A}_n$) for $(m,n)=(6,6),(8,8),(10,10),\dots$}
\end{figure}

\textbf{Lemma 2.} \emph{If $x\in\mathcal{P}(Q)$ is of spiral type 1, then }
\begin{equation}
\tilde{P}_{\pm}(x)\equiv 0,
\end{equation}
\emph{where ``$\pm$'' means that we take a doubled path with $v=\pm1$}.
\emph{This lemma can be improved if we take all irreducible subsets of the set of all arrows of $Q$, on which our polynomial vanishes.}
For instance, each minimal null-set for $4$-vertex quiver is trivial,
$6$-vertex: $\{(1,6),(3,4),(1,5),(3,6)\}$ - is not minimal, but irreducible
(recall that in this case orientation does not matter),
$8$-vertex: $\{(2,7),(3,6),(2,6),(3,5)\}$,
$12$-vertex: $\{(5,8),(3,10),(4,9),(4,8),(5,7),(3,9)\}$ and so on. 
Let's write the corresponding symmetric polynomials explicitly:
\begin{equation}
\hat{P}_{i}=\left.\frac{\tilde{P}_{+}(x_i)-\tilde{P}_{-}(x_i)}{2}\right|_{\lambda_{j,k}\mapsto \mathcal{Y}(\lambda_{j,k})}, \quad \mathcal{Y}(l_i)=\emptyset,\mathcal{Y}(\lambda>l_i)=[\lambda-l_i]
\end{equation}
\begin{equation}
\begin{aligned}
\hat{P}_4=&0
\\
\hat{P}_6=&
 \left( s_{\tiny\yng(1,1,3)_2}-s_{\tiny\yng(1,1,3)_1} \right) \Delta_{\tiny\yng(2)_3}+
 \left(s_{\tiny\yng(3)_2}- s_{\tiny\yng(3)_1}\right) \Delta_{\tiny\yng(1,2)_3}+ \left( s_
{\tiny\yng(1,3)_2}-s_{\tiny\yng(1,3)_1} \right) \Delta_{\tiny\yng(1,1,2)_3}
\\
\hat{P}_8=&
\left( s_{\tiny\yng(1,1,1,3)_3}-s_{\tiny\yng(1,1,1,3)_2} \right) \Delta_{\tiny\yng(2)_4}+
\left( s_{\tiny\yng(1,1,1,1)_3}-s_{\tiny\yng(1,1,1,3)_1} \right) \Delta_{\tiny\yng(4)_4}+
\left( -s_{\tiny\yng(3)_2}+s_{\tiny\yng(3)_3} \right) \Delta_{\tiny\yng(1,2)_4}+
\\ & +
\left( s_{\tiny\yng(1)_3}-s_{\tiny\yng(3)_1} \right) \Delta_{\tiny\yng(1,4)_4}+ \left( -s_{\tiny\yng(1,3)_2}+
s_{\tiny\yng(1,3)_3} \right) \Delta_{\tiny\yng(1,1,2)_4}+ \left( s_{\tiny\yng(1,1)_3}-s_{\tiny\yng(1,3)_1} \right) \Delta_{\tiny\yng(1,1,4)_4}+
\\ & +
\left( s_{\tiny\yng(1,1,3)_3}-s_{\tiny\yng(1,1,3)_2} \right) \Delta_{\tiny\yng(1,1,1,2)_4}+ \left( s_{\tiny\yng(1,1,1)_3}-s_{\tiny\yng(1,1,3)_1} \right) \Delta_{\tiny\yng(1,1,1,4)_4}
\\
\hat{P}_{10}=&
 \left( -s_{\tiny\yng(3)_3}+s_{\tiny\yng(3)_4} \right) \Delta_{\tiny\yng(2)_9}+ \left( -s_{\tiny\yng(5)_1}+
 s_{\tiny\yng(5)_2}+s_{\tiny\yng(5)_4}-s_{\tiny\yng(5)_3} \right) \Delta_{\tiny\yng(4)_9}+
 \\&+
 \left( s_{\tiny\yng(4)_1}-s_{\tiny\yng(2)_3} \right) \Delta_{\tiny\yng(5)_9}+ \left( -s_{\tiny\yng(1,3)_3}+
 s_{\tiny\yng(1,3)_4} \right) \Delta_{\tiny\yng(1,2)_9}+ \left( s_{\tiny\yng(1,5)_2}-s_{\tiny\yng(1,5)_3}- \right.
 \\ &
\left.-s_{\tiny\yng(1,5)_1}+s_{\tiny\yng(1,5)_4} \right)
 \Delta_{\tiny\yng(1,4)_9}+
 \left( s_{\tiny\yng(1,4)_1}-s_{\tiny\yng(1,2)_3} \right) \Delta_{\tiny\yng(1,5)_9}
+ \left( s_{\tiny\yng(1,1,3)_4}-s_{\tiny\yng(1,1,3)_3} \right) \Delta_{\tiny\yng(1,1,2)_9}+
\\ & +
  \left( -s_{\tiny\yng(1,1,5)_3}+s_{\tiny\yng(1,1,5)_2}+s_{\tiny\yng(1,1,5)_4}-s_{\tiny\yng(1,1,5)_1}
  \right) \Delta_{\tiny\yng(1,1,4)_9}+ \left( s_{\tiny\yng(1,1,4)_1}-s_{\tiny\yng(1,1,2)_3}
 \right) \Delta_{\tiny\yng(1,1,5)_9}+
\\ & + 
  \left( -s_{\tiny\yng(1,1,1,3)_3}+s_{\tiny\yng(1,1,1,3)_4}
 \right) \Delta_{\tiny\yng(1,1,1,2)_9} 
  \left( -s_{\tiny\yng(1,1,1,5)_1}+s_{\tiny\yng(1,1,1,5)_4}-
s_{\tiny\yng(1,1,1,5)_3}+s_{\tiny\yng(1,1,1,5)_2} \right) \Delta_{\tiny\yng(1,1,1,4)_9}+
\\ & +
 \left( s_{\tiny\yng(1,1,1,4)_1}-s_{\tiny\yng(1,1,1,2)_3} \right) \Delta_{\tiny\yng(1,1,1,5)_9}+
 \left( s_{\tiny\yng(1,1,1,1,3)_4}-s_{\tiny\yng(1,1,1,1,3)_3} \right) 
\Delta_{\tiny\yng(1,1,1,1,2)_9}+ 
\\ & +
\left( s_{\tiny\yng(1,1,1,1,5)_4}+s_{\tiny\yng(1,1,1,1,5)_2}-
s_{\tiny\yng(1,1,1,1,5)_3}-s_{\tiny\yng(1,1,1,1,5)_1} \right) \Delta_{\tiny\yng(1,1,1,1,4)_9}+
 \left( s_{\tiny\yng(1,1,1,1,4)_1}-s_{\tiny\yng(1,1,1,1,2)_3} \right) \Delta_{\tiny\yng(1,1,1,1,5)_9}
\\
\hat{P}_{12}= &
\left( s_{\tiny\yng(3)_5}-s_{\tiny\yng(3)_4} \right) \Delta_{\tiny\yng(2)_{11}}+
 \left( s_{\tiny\yng(5)_3}-s_{\tiny\yng(5)_4}+s_{\tiny\yng(5)_5}-s_{\tiny\yng(5)_2} \right) \Delta_{\tiny\yng(4)_{11}}+
 \left( s_{\tiny\yng(1,3)_5}-s_{\tiny\yng(1,3)_4} \right) \Delta_{\tiny\yng(1,2)_{11}}+
\\ & +
 \left( s_{\tiny\yng(1,5)_5}-s_{\tiny\yng(1,5)_2}+s_{\tiny\yng(1,5)_3}-s_{\tiny\yng(1,5)_4} \right) \Delta_{\tiny\yng(1,4)_{11}}+
 \left( -s_{\tiny\yng(1,1,3)_4}+s_{\tiny\yng(1,1,3)_5} \right) \Delta_{\tiny\yng(1,1,2)_{11}}+
\\ & +
 \left( -s_{\tiny\yng(1,1,5)_4}-s_{\tiny\yng(1,1,5)_2}+s_{\tiny\yng(1,1,5)_5}+s_{\tiny\yng(1,1,5)_3} \right) \Delta_{\tiny\yng(1,1,4)_{11}}+
 \left( -s_{\tiny\yng(1,1,1,3)_4}+s_{\tiny\yng(1,1,1,3)_5} \right) \Delta_{\tiny\yng(1,1,1,2)_{11}}+
\\ & +
 \left( -s_{\tiny\yng(1,1,1,1,5)_4}+s_{\tiny\yng(1,1,1,1,5)_5}+s_{\tiny\yng(1,1,1,1,5)_3}-s_{\tiny\yng(1,1,1,1,5)_2} \right) \Delta_{\tiny\yng(1,1,1,1,4)_{11}}+ \left( s_{\tiny\yng(1,1,1,1,1,3)_5}-s_{\tiny\yng(1,1,1,1,1,3)_4} \right) \Delta_{\tiny\yng(1,1,1,1,1,2)_{11}}+
\\ & + 
 \left( -s_{\tiny\yng(1,1,1,1,1,5)_2}-s_{\tiny\yng(1,1,1,1,1,5)_4}+s_{\tiny\yng(1,1,1,1,1,5)_5}+s_{\tiny\yng(1,1,1,1,1,5)_3} \right)
\Delta_{\tiny\yng(1,1,1,1,1,4)_{11}}+\dots
\end{aligned}
\end{equation}
Another (stand-alone) example:
\begin{equation}
-I\hat{P}_6(spiral_2)=s_{\tiny\yng(1,1,2)_2}+s_{\tiny\yng(1,2)_0}+s_{\tiny\yng(1,1,2)_4}-s_{\tiny\yng(1,1,2)_3}-s_{\tiny\yng(1,2)_5}-s_{\tiny\yng(2)_1}
\end{equation}
Let $(m,n)=(6,6)$; consider the following two paths and their concatenation:
\begin{figure}[h!]
\centering
\includegraphics[scale=1]{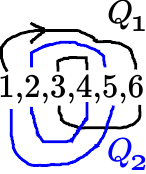}
\end{figure}
{\small{
\begin{equation}
\begin{aligned}
\tilde{P}(Q_1)= &\left( s_{\tiny\yng(1)_0}+s_{\tiny\yng(1,1)_3}+s_{\tiny\yng(1)_5}+s_{\tiny\yng(2)_2} \right) \Delta_{\tiny\yng(1)_5}+ \left( -s_{\tiny\yng(2)_0}+s_{\tiny\yng(1,2)_3}+s_{\tiny\yng(2)_5}-s_{\tiny\yng(3)_2} \right) \Delta_{\tiny\yng(2)_5}+
\\ & +
 \left( s_{\tiny\yng(1)_2}+s_{\tiny\yng(3)_0}+s_{\tiny\yng(3)_5}+s_{\tiny\yng(1,3)_3} \right) \Delta_{\tiny\yng(3)_5}+ \left( s_{\tiny\yng(1,1)_5}+s_{\tiny\yng(1,1)_0}+s_{\tiny\yng(1,1,1)_3}+s_{\tiny\yng(1,2)_2} \right) \Delta_{\tiny\yng(1,1)_5}+
\\ & + 
  \left( -s_{\tiny\yng(1,3)_2}+s_{\tiny\yng(1,1,2)_3}+s_{\tiny\yng(1,2)_5}-s_{\tiny\yng(1,2)_0} \right) \Delta_{\tiny\yng(1,2)_5}+ \left( s_{\tiny\yng(1,3)_0}+s_{\tiny\yng(1,1)_2}+s_{\tiny\yng(1,3)_5}+s_{\tiny\yng(1,1,3)_3} \right) \Delta_{\tiny\yng(1,3)_5}+
\\ & +  
   \left( s_{\tiny\yng(1,1,1)_5}+s_{\tiny\yng(1,1,2)_2}+s_{\tiny\yng(1,1,1)_0}+s_{\tiny\yng(1)_3} \right) \Delta_{\tiny\yng(1,1,1)_5}+ \left( s_{\tiny\yng(2)_3}-s_{\tiny\yng(1,1,3)_2}+s_{\tiny\yng(1,1,2)_5}-s_{\tiny\yng(1,1,2)_0} \right) \Delta_{\tiny\yng(1,1,2)_5}+
\\ & +   
    \left( s_{\tiny\yng(1,1,3)_5}+s_{\tiny\yng(3)_3}+s_{\tiny\yng(1,1,1)_2}+s_{\tiny\yng(1,1,3)_0} \right) \Delta_{\tiny\yng(1,1,3)_5}
    \quad \text{\# 5 is the longest arc}
\\
\tilde{P}(Q_2)= &
\left( s_{\tiny\yng(1)_4}+s_{\tiny\yng(1)_5}-s_{\tiny\yng(2)_2}+s_{\tiny\yng(1)_1} \right) \Delta_{\tiny\yng(1)_4}+ \left( -s_{\tiny\yng(2)_5}-s_{\tiny\yng(3)_1}+s_{\tiny\yng(3)_2}+s_{\tiny\yng(2)_4} \right) \Delta_{\tiny\yng(2)_4}+
\\ & +
 \left( -s_{\tiny\yng(3)_5}+s_{\tiny\yng(1)_2}+s_{\tiny\yng(3)_4}+s_{\tiny\yng(2)_1} \right) \Delta_{\tiny\yng(3)_4}+ \left( -s_{\tiny\yng(1,1,2)_2}+s_{\tiny\yng(1,1,1)_1}+s_{\tiny\yng(1,1)_4}+s_{\tiny\yng(1,1,1)_5} \right) \Delta_{\tiny\yng(1,1)_4}+
\\ & + 
  \left( -s_{\tiny\yng(1,1,2)_5}+s_{\tiny\yng(1,1,3)_2}+s_{\tiny\yng(1,2)_4}-s_{\tiny\yng(1,1,3)_1} \right) \Delta_{\tiny\yng(1,2)_4}+ \left( s_{\tiny\yng(1,1,1)_2}-s_{\tiny\yng(1,1,3)_5}+s_{\tiny\yng(1,1,2)_1}+s_{\tiny\yng(1,3)_4} \right) \Delta_{\tiny\yng(1,3)_4}+
\\ & +
   \left( s_{\tiny\yng(1,1)_1}+s_{\tiny\yng(1,1)_5}+s_{\tiny\yng(1,1,1)_4}-s_{\tiny\yng(1,2)_2} \right) \Delta_{\tiny\yng(1,1,1)_4}+ \left( -s_{\tiny\yng(1,3)_1}-s_{\tiny\yng(1,2)_5}+s_{\tiny\yng(1,1,2)_4}+s_{\tiny\yng(1,3)_2} \right) \Delta_{\tiny\yng(1,1,2)_4}+
\\ & +
\left( s_{\tiny\yng(1,1)_2}+s_{\tiny\yng(1,2)_1}+s_{\tiny\yng(1,1,3)_4}-s_{\tiny\yng(1,3)_5} \right) \Delta_{\tiny\yng(1,1,3)_4}
\quad \text{\# 4 is the longest arc}
\\
\tilde{P}(Q_{1,2})= &
\left( -s_{\tiny\yng(1,1,3)_2}+s_{\tiny\yng(1,1,3)_1}-s_{\tiny\yng(1,2)_4}-s_{\tiny\yng(1,1,2)_0}+s_{\tiny\yng(2)_3}+s_{\tiny\yng(1,1,2)_5} \right)
 \Delta_{\tiny\yng(2)_3}+
  \left( s_{\tiny\yng(2)_5}+s_{\tiny\yng(3)_1}-s_{\tiny\yng(3)_2}-s_{\tiny\yng(2)_4}-s_{\tiny\yng(2)_0}+s_{\tiny\yng(1,2)_3} \right) \Delta_{\tiny\yng(1,2)_3}+ 
\\ &  +  
\left( s_{\tiny\yng(1,2)_5}-s_{\tiny\yng(1,2)_0}-s_{\tiny\yng(1,1,2)_4}+s_{\tiny\yng(1,1,2)_3}+s_{\tiny\yng(1,3)_1}-s_{\tiny\yng(1,3)_2} \right) \Delta_{\tiny\yng(1,1,2)_3}=
\\ = & \frac{1}{2}\left(\mathrm{coeff}(P(Q_1),\Delta_{\tiny\yng(1,1,2)_5})-\mathrm{coeff}(P(Q_2),\Delta_{\tiny\yng(1,2)_4})\right)\Delta_{\tiny\yng(2)_3}
+ \frac{1}{2}\left(\mathrm{coeff}(P(Q_1),\Delta_{\tiny\yng(2)_5})-\mathrm{coeff}(P(Q_2),\Delta_{\tiny\yng(2)_4})\right)\Delta_{\tiny\yng(1,2)_3}+
\\ & +
\frac{1}{2}\left(\mathrm{coeff}(P(Q_1),\Delta_{\tiny\yng(1,2)_5})-\mathrm{coeff}(P(Q_2),\Delta_{\tiny\yng(1,1,2)_4})\right)\Delta_{\tiny\yng(1,1,2)_3}
= \frac{1}{2}\left(P(Q_1)-\mathrm{permut}(P(Q_2),\Delta_{*})\right)
\\
\end{aligned}
\end{equation}
}}
(This is a particular example related to the proposition 1).
\newpage
The next goal is to build infinite series like this
\begin{equation}
Z_R(S)= \frac{1}{2k}\sum\limits_{i,j}\left(\tilde{P}_R(x_i)-\tilde{P}_R(x_j)\right),
\end{equation}
where $S= \{x_1,x_2,\dots,x_i,\dots\}\subseteq \mathcal{P}(Q_{n\to\infty}) \ \text{ is a proper union of loops}$.
\begin{figure}[h!]
\centering
\includegraphics[scale=0.5]{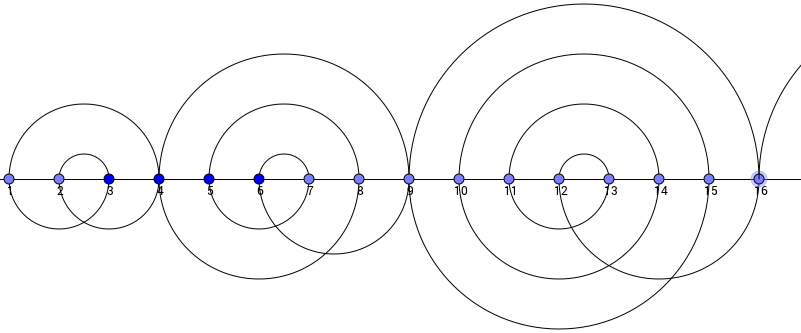}
\caption{Concatenation of $x_i\in S$ for $R=\Sigma_0$}
\end{figure}
Of course it can be rewritten in terms of Schur function expansion
\begin{equation}
Z_R(S)=\sum\limits_{\lambda}g_{\lambda}(\Delta)s_{\lambda}(p_1,p_2,\dots),
\end{equation}
or, in our \emph{``ladder notation''}, over all $\Gamma$-type \emph{weighted} Young diagrams,
using the $\mathcal{Y}$ operator
\begin{equation}
Z_R(S)=\sum\limits_{\lambda_w}g_{\lambda_w}(\Delta)\tilde{s}_{\lambda_w}(p_1,p_2,\dots), \quad w\in\mathbb{N},
\quad \lambda_w=\mathcal{Y}(\lambda)
\end{equation}
Take $S'=\{x_1,x_2\}\subset \mathcal{P}(Q_{10}), (m,n)=(10,10)$.
In this case the double-path polynomial is non-trivial (!)
\begin{equation}
\begin{aligned}
\tilde{P}_{\pm}(S')= &
\left( s_{\tiny\yng(1)_2}+s_{\tiny\yng(1,5)_8}+s_{\tiny\yng(1)_6}+s_{\tiny\yng(1)_3}+
s_{\tiny\yng(1,5)_9}+s_{\tiny\yng(1)_1}+s_{\tiny\yng(1)_7}+s_{\tiny\yng(1,1)_5}+s_{\tiny\yng(1,1)_4} \right)
 \Delta_{\tiny\yng(1)_6}+
\\ & + 
  \left( s_{\tiny\yng(1,1)_3}+s_{\tiny\yng(1,1)_2}+s_{\tiny\yng(1,1,1)_5}+
s_{\tiny\yng(1,1)_1}+s_{\tiny\yng(1,1,5)_8}+s_{\tiny\yng(1,1,1)_4}+s_{\tiny\yng(1,1,1)_6}+
s_{\tiny\yng(1,1)_7}+s_{\tiny\yng(1,1,5)_9} \right) \Delta_{\tiny\yng(1,1,1)_6}+
\\ & + 
 \left( s_{\tiny\yng(1,1,1,1)_5}+
s_{\tiny\yng(1,1,1,1)_6}+s_{\tiny\yng(1,1,1)_2}+s_{\tiny\yng(1,1,1,1)_4}+s_{\tiny\yng(1,1,1)_3}+
s_{\tiny\yng(1,1,1,5)_9}+s_{\tiny\yng(1,1,1,1,5)_8}+s_{\tiny\yng(1,1,1,1)_7}+
s_{\tiny\yng(1,1,1)_1} \right) \Delta_{\tiny\yng(1,1,1,1)_6}+
\\ & + 
 \left( s_{\tiny\yng(1,1,1,1)_3}+
s_{\tiny\yng(1,1,1,1,5)_9}+s_{\tiny\yng(1,1,1,1)_2}+s_{\tiny\yng(1,1,1,1,1)_7}+
s_{\tiny\yng(1,1,1,1,1)_5}+s_{\tiny\yng(1,1,1,1,1)_4}+s_{\tiny\yng(1,1,1,5)_8}+s_{\tiny\yng(1,1,1,1)_1}+
s_{\tiny\yng(1,1,1,1,1)_6} \right) \Delta_{\tiny\yng(1,1,1,1,1)_6}+
\\ & + 
s_{\tiny\yng(1,1,1,1,1)_3}+
s_{\tiny\yng(1,1,1,1,1)_1}+s_{\tiny\yng(1,1)_6}+s_{\tiny\yng(5)_9}+s_{\tiny\yng(1,1,1,1,1)_2}+
s_{\tiny\yng(1,1,1)_7}+s_{\tiny\yng(5)_8}+s_{\tiny\yng(1)_5}
\end{aligned}
\end{equation}
(Here we assumed $\Delta_{{\tiny{[6, 7, 8, 9, 10, 16, 17, 18, 19, 20,\dots, 96, 97, 98, 99, 100]}}} = 1$ to simplify the formula).
Its configuration is drawn on the figure \ref{fig1} (right).
\begin{figure}[h!]\label{fig1}
\centering
\includegraphics[scale=0.4]{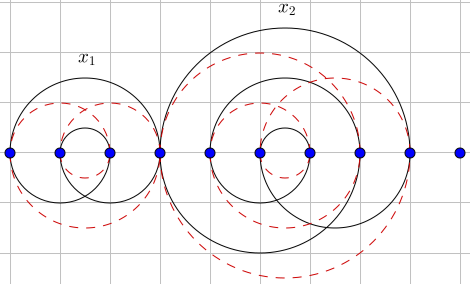}
\caption{Concatenation of $x_1,x_2\in S'$ for $R=\Sigma_0$: compatible case}
\end{figure}
Now let's extend the picture for $x_3$:
The problem is that in this case $(m,n)=(16,16)$ the image is empty, so we should add some
isolated points to the graph \ref{fig1} (left).
\newpage
Let
\begin{equation}
X=\{G_{i_1}^{(j_1)},G_{i_2}^{(j_2)},\dots,G_{i_s}^{(j_s)}\}^k \in \mathcal{P}
\end{equation}
be an element of path algebra of $Q_{\mathrm{univ}}(k)$, built as concatenation normal loops $G_{i_1},G_{i_2}$,
each one shifted on the $j$-th position. 
For example, consider
3-union $\{G_4^{(0)},G_6^{(3)},G_{8}^{(2)}\}^{10}$. QR has unique solution;
furthermore, we take $\Delta_{\tiny\yng(1,1,1)_7}=1$ (denominator) to homogenise the resulting formula.
Here is the associated quiver:
\begin{equation}
Q(X)=\begin{xy}
 0;<1pt,0pt>:<0pt,-1pt>:: 
(0,53) *+{0} ="0",
(53,53) *+{1} ="1",
(53,146) *+{2} ="2",
(53,0) *+{3} ="3",
(220,0) *+{4} ="4",
(220,146) *+{5} ="5",
(293,53) *+{6} ="6",
(147,0) *+{7} ="7",
(147,53) *+{8} ="8",
(147,146) *+{9} ="9",
"2", {\ar"0"},
"0", {\ar"3"},
"1", {\ar"2"},
"3", {\ar"1"},
"8", {\ar"2"},
"2", {\ar"9"},
"7", {\ar"3"},
"3", {\ar"8"},
"6", {\ar"4"},
"4", {\ar"7"},
"8", {\ar"4"},
"5", {\ar"6"},
"8", {\ar"5"},
"9", {\ar"5"},
\end{xy}
\end{equation}
\begin{figure}[h!]\label{fig1}
\centering
\includegraphics[scale=0.5]{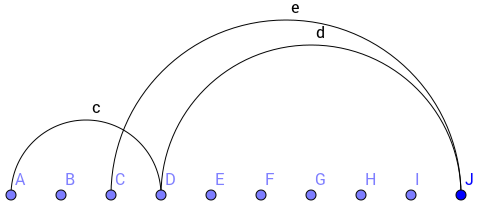} \\
\includegraphics[scale=0.5]{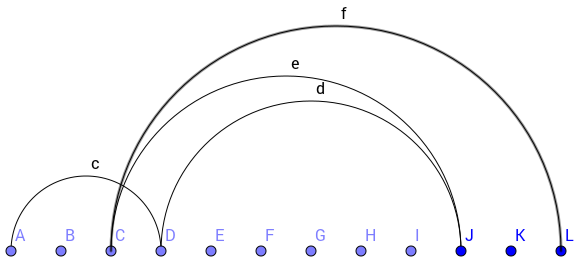}
\caption{Non-trivial concatenation of loops $X=\{G_4^{(0)},G_6^{(3)},G_{8}^{(2)}\}^{10}$ and
$X'=\{G_4^{(0)},G_6^{(3)},G_{8}^{(2)},\underline{G_{10}^{(2)}}\}^{10}$ (only the largest arcs are shown)}
\end{figure}
Then $\tilde{P}_{\pm}(X)$ become homogeneous in $\Delta$:
\begin{equation}
\begin{aligned}
\tilde{P}_{\pm}(X)=& \left( s_{\tiny\yng(1,5)_{9}}+s_{\tiny\yng(1)_{0}}+s_{\tiny\yng(1,1)_{5}}+s_{\tiny\yng(1)_{1}}+s_{\tiny\yng(1)_{2}}+
s_{\tiny\yng(1,1)_{4}}+s_{\tiny\yng(1)_{6}}+s_{\tiny\yng(1)_{3}}+s_{\tiny\yng(1,5)_{8}}+s_{\tiny\yng(1)_{7}} \right) 
\Delta_{\tiny\yng(1)_{6}}+
\\ &
+ \left( s_{\tiny\yng(1,1)_{0}}+s_{\tiny\yng(1,1,5)_{8}}+s_{\tiny\yng(1,1)_{3}}+
s_{\tiny\yng(1,1)_{1}}+s_{\tiny\yng(1,1,1)_{6}}+s_{\tiny\yng(1,1,1)_{5}}+s_{\tiny\yng(1,1)_{2}}+s_{\tiny\yng(1,1,1)_{4}}+
s_{\tiny\yng(1,1)_{7}}+s_{\tiny\yng(1,1,5)_{9}} \right) \Delta_{\tiny\yng(1,1,1)_{6}}+
\\ &
+ \left( s_{\tiny\yng(1,1,1)_{2}}+
s_{\tiny\yng(1,1,1,1)_{7}}+s_{\tiny\yng(1,1,1,1)_{6}}+s_{\tiny\yng(1,1,1)_{3}}+s_{\tiny\yng(1,1,1,1)_{5}}+
s_{\tiny\yng(1,1,1,1)_{4}}+s_{\tiny\yng(1,1,1,5)_{9}}+s_{\tiny\yng(1,1,1)_{0}}+s_{\tiny\yng(1,1,1)_{1}}+
s_{\tiny\yng(1,1,1,1,5)_{8}} \right) \Delta_{\tiny\yng(1,1,1,1)_{6}}+
\\ &
+ \left( s_{\tiny\yng(1,1,1,1,5)_{9}}+
s_{\tiny\yng(1,1,1,5)_{8}}+s_{\tiny\yng(1,1,1,1)_{1}}+s_{\tiny\yng(1,1,1,1,1)_{7}}+s_{\tiny\yng(1,1,1,1,1)_{5}}+
s_{\tiny\yng(1,1,1,1,1)_{6}}+s_{\tiny\yng(1,1,1,1)_{0}}+s_{\tiny\yng(1,1,1,1)_{2}}+s_{\tiny\yng(1,1,1,1)_{3}}+
s_{\tiny\yng(1,1,1,1,1)_{4}} \right) \Delta_{\tiny\yng(1,1,1,1,1)_{6}}+
\\ &
+ s_{\tiny\yng(1,1,1,1,1)_{0}}+
s_{\tiny\yng(1,1,1,1,1)_{1}}+s_{\tiny\yng(5)_{8}}+s_{\tiny\yng(1)_{5}}+s_{\tiny\yng(1,1,1,1,1)_{3}}+s_{\tiny\yng(5)_{9}}+
s_{\tiny\yng(1)_{4}}+s_{\tiny\yng(1,1)_{6}}+s_{\tiny\yng(1,1,1,1,1)_{2}}
\end{aligned}
\end{equation}
Remark: we can see \underline{only antisymmetric coefficients} $\Delta$ here! To extend this case for $G_{10}$
one should exclude $G_{8}^{(2)}$.

The resulting picture is that \underline{each loop can be scaled separately}, such that increasing $k,s$ simultaneously
does not change the consistency of QR. Our next aim is to investigate asymptotical properties, derived from
$\{G_{i_1}^{(j_1)},G_{i_2}^{(j_2)},\dots,G_{i_s}^{(j_s)}\}^k$ when $k,s\to\infty$.
This is of further development.
%

\end{document}